\documentclass[fleqn]{mat01}
\usepackage{times,mathtimy,amssymb,epsfig,latexsym}
\begin{document}

\setcounter{page}{141}
\firstpage{141}

\renewcommand{\theequation}{\thesection\arabic{equation}}

\def\d{\mbox{\rm d}}
\def\e{\mbox{\rm e}}

\font\zz=msam10 at 10pt
\def\Box{\mbox{\zz{\char'244}}}

\def\note{\trivlist \item[\hskip \labelsep{\it Note:}]}

\newcommand{\field}[1]{\ensuremath{\mathbb{#1}}}
\newcommand{\CC}{\field{C}}
\newcommand{\HH}{\field{H}}
\newcommand{\RR}{\field{R}}
\newcommand{\ZZ}{\field{Z}}

\newcommand{\del}{\partial}
\newcommand{\delb}{\bar\partial}

\title{A complete conformal metric of preassigned negative Gaussian
curvature for a punctured hyperbolic Riemann surface}

\markboth{Rukmini Dey}{Preassigned negative Gaussian curvature}

\author{RUKMINI DEY}

\address{Department of Mathematics, Indian Institute of Technology,
Kanpur~208~016, India\\
\noindent E-mail: rukmini@iitk.ac.in}

\volume{114}

\mon{May}

\parts{2}

\Date{MS received 1 December 2003; revised 14 February 2004}

\begin{abstract}
Let $h$ be a complete metric of Gaussian curvature $K_0$ on a punctured
Riemann surface of genus $g \geq 1$ (or the sphere with at least three
punctures). Given a smooth negative function $K$ with $K=K_0$ in
neighbourhoods of the punctures we prove that there exists a metric
conformal to $h$ which attains this function as its Gaussian curvature
for the punctured Riemann surface. We do so by minimizing an appropriate
functional using elementary analysis.
\end{abstract}


\maketitle

\section{Introduction}

Let $\Sigma$ be a punctured hyperbolic Riemann surface, namely a Riemann
surface of genus $g \geq 1$ with at least one puncture (or the sphere
with at least three punctures). Let $\d s^2 = h\, \d z \otimes \d
\bar{z} $ be a complete metric on $\Sigma$ with $K_0= K_0(z, \bar{z})$
its Gaussian curvature, negative near the punctures. Example of an
initial metric is given in \S4.

In this paper we prove the result that given an {\it arbitrary} smooth
negative function $K$, with $K=K_0$ in neighbourhoods of the punctures,
there exists a metric conformal to $h$ which attains this function as
its Gaussian curvature\footnote{Note that there are minimal surfaces
with ends, topologically equivalent to punctured Riemann surfaces, which
have negative curvature everywhere \cite{JM,O,X}.}. Our proof is
elementary, using Hodge theory, i.e., the existence of the Green's
operator for the Laplacian on a compact manifold with boundary. This
proof is a generalization of the proof given in an earlier paper that
for a compact Riemann surface of genus $g>1$ any negative function $K$
is attained as the Gaussian curvature of a metric conformal to the given
one~\cite{D}.

In~\cite{HT} Appendix A, there is a theorem.

\begin{theorem}[\!]\hskip -.5pc {\bf (Hulin--Troyanov).}\ \
Let $S^{\prime}$ be a non-compact Riemann surface of finite type. Assume
$S^{\prime}$ is neither conformally equivalent to ${\CC}$ nor
${\CC}^{*}$. Let $K\hbox{\rm :}\ S^{\prime} \rightarrow {\RR}$ be any bounded
locally H\"{o}lder-continuous function. Then{\rm ,} there exists a conformal
metric $g$ on $S^{\prime}$ with curvature $K$.
\end{theorem}

However, such metrics are not usually complete. If we restrict ourselves
to complete metrics but forgo having an example in each conformal class,
there is a result by Burago (proof is to be found in \cite{HT})
generalizing Kazdan and Warner that any function negative\newpage \noindent on an open set
with non-positive value somewhere near the punctures is achieved as a
prescribed Gaussian curvature \cite{Bu,HT,KW2}.

\begin{theorem}[\!]\hskip -.4pc {\bf (Burago--Kazhdan--Warner).}\ \
Let $K$ be a smooth function on an open surface $S^{\prime}$ of finite
type. The following conditions are necessary and sufficient for $K$
to be the curvature of some complete metric of $S^{\prime}${\rm :}

\begin{enumerate}
\renewcommand{\labelenumi}{\rm (\roman{enumi})}
\leftskip .1pc
\item $\lim \inf_{x \rightarrow \infty}\ K(x) \leq 0$ at each end of
$S^{\prime}${\rm ;}

\item when Euler characteristic of $S^{\prime} <0${\rm ,} assume also{\rm ,}
\hbox{$\inf K <0$}{\rm ;} when Euler characteristic of $S^{\prime} = 0${\rm ,}
assume also \hbox{$\inf K <0$} or $K \equiv 0$.
\end{enumerate}
\end{theorem}

The proof of this theorem is also given in \cite{HT}. However they do
not find a metric in every conformal class which has this property.

Our proof is very different and our result is new in the sense that we
show that if the initial metric has negative Gaussian curvature near the
punctures and given any negative function with the same value as the
initial one near the punctures, there exists a metric in the conformal
class of the initial metric which has this function as its Gaussian
curvature. Moreover, the behaviour of the initial metric near the
punctures do not change (the conformal equivalence factor
$\hbox{e}^{\sigma}$ is 1 near the punctures). Thus if we started with a
complete metric, completeness is preserved. If the {\it initial} metric
has Gaussian curvature the function $K_0$, negative somewhere, such
that $\lim_{x \rightarrow \infty} K_0 = f_0 <0 $ (i.e. its value is
$f_0$ near the punctures) then by Burago--Kazhdan--Warner's result we
would have a complete metric which achieves $K_0$ as its Gaussian
curvature. Let $K <0$ be another function such that $\lim_{x \rightarrow
\infty} K = f_0 $ then our result shows that there exists a metric
conformal to the previous one such that $K$ is achieved as its Gaussian
curvature. This indicates that there could be many conformal classes
which contain a metric with $K$ as its Gaussian curvature. In the last
section we give an example of the initial metric which has negative
Gaussian curvature near the punctures.

By Thereom 2.1a, we also have the result that any arbitrary negative
function $K$ is allowed as the Gaussian curvature of an unpunctured
Riemann surface of genus $g \geq 1$ with at least one boundary component
or a sphere with at least three boundary components.

For prescribed curvature on surfaces with conical singularities, see
~\cite{T}.

Let $p_i$, $i=1,\ldots,n$ be the punctures on $\Sigma$. Let us choose a
disc $D_i$ about the point $p_i$ with coordinates $ (r_i, \theta_i) $
such that $D_i = \{ (r_i, \theta_i)$\,: $0< r_i < 1,\ 0< \theta_i \leq 2
\pi \}$ with the point $p_i$ corresponds to $r_i=0$. No other puncture
lies in this neighbourhood. Let $M = \Sigma - \cup_{i=1}^n D_i$. Thus
$M$ is a compact manifold with boundary $\partial M = \cup_{i=1}^n
\partial D_i$. We assume $K_0$ is negative on $\partial M$. In \S4 we
show by an example that this is not a stringent condition on the metric
$h$. Also $M$ is of finite volume with respect to the metric $h$ and
$K_0$ is smooth on $M$. $K_0$ may have singularities at the punctures.
Let $K <0 $ be a given function on $M$ and $K = K_0 $ on $\cup_{i=1}^n
D_i$. On $M$ we show that there exists a metric conformal to $h$ which
has $K$ as its Gaussian curvature. We show this by minimizing the
functional
\begin{equation*}
S[\sigma] = \int_M (K(\sigma)-K)^2 \e^{2\sigma}\d\mu
\end{equation*}
over functions in class $T$ where $T= \{ \sigma \in
C^{2}(M,\RR)|\sigma|_{\partial M} =\partial_{\nu} \sigma|_{\partial M} =
\Delta \sigma|_{\partial M} =0 \}^{-}$ where $-$ denotes the closure in
$W^{2,2}$. Here $K(\sigma)$ stands for the Gaussian curvature of the
metric $\e^{\sigma} h$, and $\d\mu = \frac{\sqrt{-1}}{2} h\ \d z
\wedge\d\bar{z}$ is the area form for the metric $h$. Using the Sobolev
embedding theorem for compact manifolds with boundary we show that
$S[\sigma]$ takes its absolute minimum, zero, on $C^{\infty}(M)$ which
corresponds to a metric on $M$ of negative curvature $K$. Then we extend
the minimizer $\sigma$ by zero on $D_i$ and thus near the punctures the
metric $\e^{\sigma} h\, \d z \otimes \d \bar{z} $ and the curvature $K$
remain $ h\, \d z \otimes \d \bar{z}$ and $K_0$ respectively.

\section{The main theorem}

\subsection{\it All notations are as in {\rm \S1}}

The functional $S[\sigma] = \int_M
(K(\sigma) - K)^2 \e^{2\sigma} \d \mu$ is non-negative on $T$, so that
its infimum
\begin{equation*}
S_0 = \inf\{S[\sigma],\,\sigma \in T\}
\end{equation*}
exists and is non-negative. Let $\{\sigma_n\}_{n=1}^{\infty}\subset T$
be a corresponding minimizing sequence,
\begin{equation*}
\lim\limits_{n\rightarrow\infty} S[\sigma_n] =S_0.
\end{equation*}
Our main result is the following\vspace{.2pc}

\begin{theorem}[\!]
{\rm (a)} Let $M = \Sigma - \cup_{i=1}^n D_i$ be a compact Riemann
surface with boundary{\rm ,} with at least three boundary components if
genus $g=0$. Let $h\, \d z \wedge \d \bar{z}$ be its initial metric and
$K_0$ be its initial curvature such that $K_0 <0$ on $\cup_{i=1}^n D_i$.
Let $K <0$ be an arbitrary negative function with $K = K_0$ on
$\cup_{i=1}^n D_i$ {\rm (}i.e.{\rm ,} the curvature is left intact near the
punctures{\rm )}. The infimum $S_0 =0$ is attained at a unique
$\sigma\in C^{\infty}(M,\RR)${\rm ,} i.e.{\rm ,} the minimizing sequence
$\{\sigma_n\}$ contains a subsequence that converges in strong $W^{2,2}$
to a unique $\sigma\in C^{\infty} (M,\RR)$ such that $S[\sigma]=S_0 =
0$. The corresponding metric $\e^{\sigma}h\,\d z\otimes\d\bar{z}$ is a
metric on $M$ of preassigned negative curvature $K$. {\rm (b)} On
$\cup_{i=1}^n D_i${\rm ,} $\sigma$ is extended as\break zero.\vspace{.3pc}
\end{theorem}

\subsection{\it Uniform bounds}

Since $\{\sigma_n\}$ is a minimizing sequence, we have the obvious
inequality
\begin{equation}\label{1}
S[\sigma_n]=\int_{M}(K_n-K)^2\e^{2\sigma_n}\d\mu =\int_M \left(K_0 -
\frac{1}{2}\Delta_h\sigma_n -K\e^{\sigma_n}\right)^2 \d\mu\leq m
\end{equation}
for some $m>0$, where we denoted by $K_n$ the Gaussian curvature
$K(\sigma_n)$ of the metric $\e^{\sigma_n}h$ and by $K_0$ that of the
metric $h$, and used that
\begin{equation*}
K_n=\e^{-\sigma_n}\left(K_0 - \frac{1}{2} \Delta_h \sigma_n\right).
\end{equation*}
\begin{note}
Here $\Delta_h = 4h^{-1} ({\del^2}/{\del z\del\bar{z}})$
stands for the Laplacian defined by the metric $h$ on $M$.\vspace{.2pc}
\end{note}

\setcounter{defin}{1}
\begin{lemma}
There exists a constant $C_1$  such that{\rm ,} uniformly in $n${\rm ,}
\begin{equation*}
\int_M (\Delta_h\sigma_n)^2 \d\mu <C_1.
\end{equation*}
\end{lemma}
\vspace{.3pc}

\begin{proof}
By Minkowski inequality, and using \eqref{1}, we get
\begin{align*}
\left[\! \int_M \left(\! -\frac{1}{2} \Delta_h\sigma_n - K\e^{\sigma_n}\!\right)^2 \d \mu\!\right]^{1/2}\! &\leq
\left[\int_M \left(K_0 - \frac{1}{2}\Delta_h\sigma_n -
K\e^{\sigma_n}\right)^2 \d \mu\right]^{1/2}\\
&\quad\, + \left[\int_M (K_0)^2 \d \mu\right ]^{1/2}\leq m^{1/2}+c=C,
\end{align*}
so that
\begin{equation}\label{2}
\frac{1}{4}\int_M (\Delta_h\sigma_n)^2 \d\mu+ \int_M K^2 \e^{2\sigma_n}\d\mu +
\int_M \Delta_h\sigma_n \e^{\sigma_n} K \d\mu \leq C^2.
\end{equation}

We will show that
\begin{equation}
\int_M K^2 \e^{2\sigma_n}\d\mu + \int_M \Delta_h\sigma_n \e^{\sigma_n} K \d \mu
= B_{n1} + B_{n2} + B_{n3},
\end{equation}
where $B_{n1} \geq 0$, $B_{n2} \geq 0$, $|B_{n3}| \leq 3 D^2$ where $D^2$ is
a constant independent of $n$.

From (2.3) and (2.2) the result  will follow since we will have
\begin{align*}
C^2 + 3D^2 &\geq C^2 - B_{n3} \geq \frac{1}{4}\int_M (\Delta_h\sigma_n)^2 \d\mu
+ B_{n1} + B_{n2}\\[.2pc]
&\geq \frac{1}{4}\int_M (\Delta_h\sigma_n)^2 \d\mu.
\end{align*}
Just renaming the constants, we will have the result.

Integrating by parts we get
\begin{align*}
 \int_M \Delta_h\sigma_n \e^{\sigma_n} K \d \mu &= - \int_M
|\partial_{z} \sigma_n|^2 \e^{\sigma_n} K \d \mu -  \int_M (\partial_z
\sigma_n) (\partial_{\bar{z}} K) \e^{\sigma_n} \d \mu\\[.2pc]
&\quad\, -\int_{\partial M} (\partial_{\nu} \sigma_n) K \e^{\sigma_n} \d \mu\\[.2pc]
&=  \int_M |\partial_{z} \sigma_n|^2 \e^{\sigma_n} |K|\d \mu -  \int_M (\partial_z \sigma_n) g |K|\ \e^{\sigma_n} \d \mu
\end{align*}
since $K$ is negative, $\partial_{\nu} \sigma_{n} |_{\partial M} = 0$ and
where we define $g= {\partial_{\bar{z}} K}/{|K|}$.\vspace{.2pc}

Let $M = \Omega_{n1} \cup \Omega_{n2} \cup \Omega_{n3}$, a disjoint union of
sets defined as follows:\vspace{.2pc}

On $\Omega_{n1}$, (1) $|\partial_{z} \sigma_n |> |g|$.\vspace{.4pc}

On $\Omega_{n2}$, (2) $|\partial_{z} \sigma_n |\leq |g|$ and
$ |K|\ \e^{\sigma_n} > |g|^2$.\vspace{.4pc}

On $\Omega_{n3}$, (3) $|\partial_{z} \sigma_n |\leq |g|$ and
$ |K|\ \e^{\sigma_n} \leq |g|^2$.\vspace{.4pc}

Let $B_{ni} = \int_{\Omega_{ni}} K^2 \e^{2\sigma_n}\d\mu + \int_{\Omega_{ni}}
\Delta_h\sigma_n \e^{\sigma_n} K \d \mu,\ i=1,2,3.$\vspace{.3pc}

We will show that $B_{n1} \geq 0$.
\begin{align*}
B_{n1} &= \int_{\Omega_{n1}} K^2 \e^{2\sigma_n} \d \mu + \int_{\Omega_{n1}} |\partial_{z} \sigma_n|^2 |K|\ \e^{\sigma_n} \d \mu
- \int_{\Omega_{n1}} (\partial_{z} \sigma_n) g |K|\ \e^{\sigma_n} \d \mu\\[.2pc]
&\geq  \int_{\Omega_{n1}} K^2 \e^{2\sigma_n} \d \mu + \int_{\Omega_{n1}} |\partial_{z} \sigma_n|^2 |K|\ \e^{\sigma_n} \d \mu
- \int_{\Omega_{n1}} |\partial_{z} \sigma_n| |g| |K|\ \e^{\sigma_n} \d \mu%\\[.2pc]
\end{align*}
\begin{align*}
&=\int_{\Omega_{n1}} K^2 \e^{2\sigma_n} \d \mu + \int_{\Omega_{n1}} |\partial_{z} \sigma_n| |K|\ \e^{\sigma_n}(|\partial_z \sigma_n| - |g|) \d\mu\\[.2pc]
&\geq 0
\end{align*}
by (1) in the definition of $\Omega_{n1}$.

Next we shall show that $B_{n2} \geq 0$.
\begin{align*}
B_{n2} &= \int_{\Omega_{n2}} K^2 \e^{2\sigma_n} \d \mu + \int_{\Omega_{n2}} |\partial_{z} \sigma_n|^2 |K|\ \e^{\sigma_n} \d \mu - \int_{\Omega_{n2}} (\partial_{z} \sigma_n) g |K|\ \e^{\sigma_n} \d \mu\\[.2pc]
&\geq  \int_{\Omega_{n2}} K^2 \e^{2\sigma_n} \d \mu + \int_{\Omega_{n2}} |\partial_{z} \sigma_n|^2 |K|\ \e^{\sigma_n} \d \mu - \int_{\Omega_{n2}} |\partial_{z} \sigma_n| |g| |K|\ \e^{\sigma_n} \d \mu\\[.2pc]
&=\int_{\Omega_{n2}} |K|\ \e^{\sigma_n} (|K|\ \e^{\sigma_n} - |\partial_{z} \sigma_n| |g|) \d \mu + \int_{\Omega_{n2}} |\partial_z \sigma_n|^2 \e^{\sigma_n} |K| \d \mu\\[.2pc]
&\geq \int_{\Omega_{n2}} |K|\ \e^{\sigma_n} (|K|\ \e^{\sigma_n} - |g|^2) \d \mu +
\int_{\Omega_{n2}} |\partial_z \sigma_n|^2 \e^{\sigma_n} |K| \d \mu\\[.2pc]
&\geq 0,
\end{align*}
by using the two conditions (2) defining $\Omega_{n2}$.

Next we shall show that $B_{n3}$ is uniformly bounded.
\begin{align*}
|B_{n3}| &\leq \int_{\Omega_{n3}} K^2 \e^{2\sigma_n} \d \mu +
\int_{\Omega_{n3}} |\partial_{z} \sigma_n|^2 |K|\ \e^{\sigma_n} \d \mu +
\int_{\Omega_{n3}} |\partial_{z} \sigma_n| |g| |K|\ \e^{\sigma_n} \d
\mu\\[.2pc]
&\leq 3 D^2,
\end{align*}
where $D^2 = \max |g|^4 \mu(M)$, where $\mu(M)$ is the volume of $M$.
This follows from the two conditions $|K|\ \e^{\sigma_n} \leq |g|^2$ and
$ |\partial_{z} \sigma_n| \leq |g|$ on $\Omega_{n3}$. $D^2$ is a finite
constant ($\max |g|^4$ is finite since $K$ is non-zero and the volume of
$M$ is finite), independent of $n$. Thus the result follows.\hfill$\Box$\vspace{.3pc}
\end{proof}

\subsection{\it Pointwise convergence of $\sigma_n$}

\setcounter{defin}{2}
\begin{lemma}
$\{ \sigma_n \}$ is uniformly bounded in $W^{2,2}(M)$.
\end{lemma}

\begin{proof}
Let us recall that there is a Green's function  on the compact manifold
$M$ with boundary such that if $u \in C^{2}(M)$ then
\begin{equation*}
-u(x) = \int \int_{M} G(x,z) \Delta u(z) \d \mu(z) +
\int_{\partial M} \frac{\partial G}{\partial \nu}(x, w) u (w) \d l
\end{equation*}
(see for e.g.~\cite{C}, p.~174).  In particular if $u|_{\partial M} = 0$
then
\begin{equation*}
-u(x)= \int \int_{M} G(x,z) \Delta u (z) \d \mu(z).
\end{equation*}
Thus there exists a Green's operator $G$ such that $\Delta G u = -u$ and
$G \Delta u = - u$ (see ~\cite{C}, p.~177). Let $u = \sigma_n.$ Since
$\sigma_n |_{\partial M} =0$, $G \Delta \sigma_n = -\sigma_n$.

Next we will show that since by Lemma~2.2, $\Delta \sigma_n$ is uniformly
bounded in $L^2(M)$, $\sigma_n$ is uniformly bounded in $W^{2,2}(M)$.

Let $\tau_n = \Delta \sigma_n$. We know that $\tau_n|_{\partial M}=0$
and is uniformly bounded in $L^2(M)$, i.e., $\| \tau_n \|_2 < C, \forall
n$. Then $\tau_n$ can be written as $\tau_n = \sum_{i=1}^{\infty} c_n^i
\phi_i $ where $\phi_i$ is a solution to the equation $\Delta_h \phi +
\lambda_i \phi = 0$, with $\phi|_{\partial M} = 0$, (\cite{C},
pp.~8--9). $\{\phi_i\}_{i=1}^{\infty}$ is an orthonormal sequence of
eigenfunctions of the Laplacian. In fact since we are in the Dirichlet
case, the least eigenvalue $\lambda_1 >0$. We will show that since $\{
\tau_n \} $ is uniformly bounded in $L^2(M)$, $\{ G\tau_n \}$ is
uniformly bounded in $W^{2,2}(M)$. $G \tau_n = \sum_{i=1}^{\infty} c_n^i
G \phi_i$. Since $\phi_i|_{\partial M} = 0$, $G \Delta \phi_i = -
\phi_i$. Since $\Delta_h \phi_i = -\lambda_i \phi_i$, $G \phi_i =
{\phi_i}/{\lambda_i}$. Thus $G\tau_n = \sum_{i=1}^{\infty}
({c_n^i}/{\lambda_i}) \phi_i$. Since $\lambda_i \rightarrow \infty$ and
since $\| \tau_n \|_2 = (\sum_{i=1}^{\infty} c_n^{i2})^{1/2} < C$, $\| G
\tau_n \|_2 = (\sum_{i=1}^{\infty} ({c_n^i}/{\lambda_i})^2)^{1/2} < C,$
$\forall n$ if $\lambda_1 >1$. If $\lambda_1 <1$, $\| G \tau_n \|^2_2
\leq ({C^2}/{\lambda_1^2})$ $\forall n$. Thus $\|G \tau_n\| \leq C_2\
\forall n$ where
\begin{align*}
C_2 &= C\quad \hbox{if}\ \ \lambda_1 >1\\
&= \frac{C}{\lambda_1}\quad \hbox{if}\ \ \lambda_1 <1.
\end{align*}

Next we show that the first and second derivative norms of $G \tau_n$
are uniformly bounded in $L^2(M)$. The second derivative is defined as
$\nabla ^2 \tau = \nabla^1 \nabla \tau$ where $\nabla \tau $ is the
usual covariant derivative of $f$ taking values in $T^*M$ and $\nabla^1$
is the covariant derivative of sections of $T^*M$.

First we show that the first derivative norm is uniformly bounded. For
this we observe that the norm of the Laplacian is bounded uniformly.
This is because $\Delta_h G \tau_n = -\tau_n $ and $\tau_n$ is uniformly
bounded in $L^2(M)$. Thus $\| \Delta_h G \tau_n \|_2 < C$.

To show that  $\| \nabla G \tau_n \|^2_2 <C^2$ we note that since $G\tau_n$
is zero on $\partial M$,
\begin{align*}
\int_M | \nabla G \tau_n |^2 \d \mu &= - \int_M G\tau_n \Delta_h G \tau_n\\[.2pc]
& \leq \| G \tau_n \|_2 \| \Delta G \tau_n \|_2\\[.2pc]
& \leq C_2 C.
\end{align*}

Next we show that the second derivative norm $\| \nabla^2 G \tau_n \|_2
< C_1$, $\forall n$. The norm of the Laplacian is bounded in $L^2(M)$.
Using Weitzenbock formula we show that then the second derivative norm
is uniformly bounded. By Weitzenbock formula ~\cite{B}, $\Delta_h =
\nabla^* \nabla + \tau(K_0) $ where $\nabla^*$ is the adjoint of the
covariant derivative. Note that $\nabla^* \nabla = - \nabla^{1} \nabla =
- \nabla^2$ (\cite{B}, p.~52). $\tau(K_0)$ is some continuous function
of the curvature $K_0$. In other words, $\Delta_h (G\tau_n) = - \nabla^2
(G\tau_n) + \tau(K_0) (G \tau_n)$. Thus by Minkowski's inequality,
\begin{align*}
\| \nabla^2 (G \tau_n) \|_2 &\leq \| \Delta_h (G \tau_n) \|_2 +
\| \tau(K_0) G \tau_n \|_2\\[.2pc]
& \leq C + \min_M |\tau(K_0)| \| G \tau_n \|_2\\[.2pc]
& \leq C + \min_M |\tau(K_0)| C_2\\[.2pc]
& \leq C_1.
\end{align*}
Since $M$ is a compact manifold with boundary, $\min_M |\tau(K_0)|$ is bounded.

Thus $G\tau_n $ is uniformly bounded in $W^{2,2}$ since $\| G \tau_n\|_2$,
$\| \nabla (G \tau_n)\|_2$, $\| \nabla^2 (G \tau_n)\|_2 $ are uniformly bounded.
Since $\tau_n = \Delta_h \sigma_n$ and $G \Delta_h \sigma_n
= \sigma_n$, $\sigma_n$ are uniformly bounded in $W^{2,2}$.\hfill $\Box$\vspace{.4pc}
\end{proof}

Now we can formulate the main result of this subsection.

\newpage

\begin{proposition}$\left.\right.$\vspace{.5pc}

\noindent The sequence $ \{\sigma_n \}_{n=1}^{\infty}$
contains a subsequence $ \{ \sigma_{l_n}\}_{n=1}^{\infty}$ with
the following properties{\rm :}

\begin{enumerate}
\renewcommand{\labelenumi}{\rm (\alph{enumi})}
\item The sequences $\{ \sigma_{l_n} \}_{n=1}^{\infty}$ and
$\{\e^{\sigma_{l_n}}\}$ converge in $C^0(M)$ topology to continuous
functions $\sigma$ and $\e^{\sigma}$ respectively. Moreover{\rm ,}
$\sigma \in W^{2,2}(M)$.

\item The subsequence $\{\Delta_h \sigma_{l_n} \}$ converges
weakly in $L^2$ to $f\doteq \Delta_h^{\rm distr} \tilde{\sigma}$ -- a
distribution Laplacian of $\sigma$.

\item Passing to this subsequence $\{\sigma_{l_n} \}${\rm ,} the
following limits exist{\rm :}
\begin{align*}
\hskip -.55cm &\lim\limits_{n\rightarrow\infty}\|\Delta_h
\sigma_{l_n}\|_2 =\| \Delta_h^{\rm distr} \sigma \|_2,\\
\hskip -.55cm &\lim\limits_{n\rightarrow\infty}S[\sigma_{l_n}]
=S_0=\int_M \left( K_0-\frac{1}{2} \Delta_h^{\rm distr} \sigma -K
\e^{\sigma} \right)^2 \d\mu.
\end{align*}
In fact{\rm ,} the convergence in {\rm (b)} is strong in $L^2$.
\end{enumerate}
\end{proposition}

\begin{proof}
Part (a) follows from the Sobolev embedding theorem and Rellich
lemma for compact manifolds with boundary, since, for $\dim M =
2$, the space $W^{2,2}(M)$ is compactly embedded into $C^0(M)$
(see, e.g.~\cite{A}). Therefore the sequence $\{\sigma_n\}$,
which  according to Lemma~2.3, is uniformly bounded in
$W^{2,2}(M)$, contains a convergent subsequence in $C^0(M)$.
Passing to this subsequence $\{ \sigma_{l_n} \}$ we can assume
that there exists a function $\sigma \in C^0(M)$ such that
\begin{equation*}
\lim\limits_{n\rightarrow\infty}\sigma_{l_n}=\sigma.
\end{equation*}
Since $\sigma_n$'s are uniformly bounded in a Hilbert space
$W^{2,2} (M)$, they weakly converge to $s\in W^{2,2}(M)$ (after
passing to a subsequence if necessary). The uniform limit
coincides with $s$  so that $\sigma=s\in W^{2,2}(M)$.
\end{proof}

In order to prove (b), set $\psi_n=\Delta_h \sigma_{l_n}$ and
observe that, according to part (a) of Lemma~2.2, the sequence
$\{\psi_n\}$ is bounded in $L^2$. Therefore, passing to a
subsequence, if necessary, there exists $f\in L^2(M)$ such that
\begin{equation*}
\lim\limits_{n\rightarrow\infty}\int_M\psi_n g=\int_M fg
\end{equation*}
for all $g\in L^2(M)$. In particular, considering $g\in
C^{\infty}(M)$, this implies $f =\Delta_h^{\rm distr} \sigma$.

In order to prove (c) we use the following lemma.\vspace{.2pc}

\begin{lemma}
If a sequence $\{\psi_n\}$ converges to $f\in L^2$ in the weak
topology{\rm ,} then
\begin{equation*}
\lim\limits_{n\rightarrow\infty}\|\psi_n\|\geq\|f\|.
\end{equation*}
Further $\lim_{n \rightarrow \infty} \| \psi_n \| = \|f \| $ iff there is
strong convergence.\vspace{.2pc}
\end{lemma}

\begin{proof}
The lemma follows from considering the following inequality:
\begin{equation*}
\lim\limits_{n \rightarrow \infty} \int ( \psi_n - f )^2 \d \mu \geq 0.
\end{equation*}

$\left.\right.$\vspace{-2.6pc}

\hfill$\Box$\vspace{1.4pc}

To continue with the proof of the proposition, suppose $\lim_{n\rightarrow
\infty}\|\psi_n\|>\|f\|$. Using the definition of the functional, we have
\begin{align*}
S[\sigma_n] &= \int_M \left(K_0 - \frac{1}{2}\Delta_h \sigma_n
-K\e^{\sigma_n}\right)^2 \d\mu\\[.2pc]
&= \frac{1}{4}\|\psi_n\|^2+\|K_0
-K\e^{\sigma_n}\|^2 - \int_M \psi_n(K_0-K\e^{\sigma_n})\d\mu.
\end{align*}
From parts (a) and (b) it follows that the sequence $S[\sigma_n]$
converges to $S_0$ and
\begin{align*}
S_0 &=\lim\limits_{n\rightarrow\infty} S[\sigma_n]\\[.2pc]
&=\lim\limits_{n\rightarrow \infty}\frac{1}{4}\|\psi_n\|^2+\|K_0-K\e^{\sigma}\|^2-
 \int_M f(K_0-K\e^{\sigma})\d\mu\\[.2pc]
&\geq\frac{1}{4}\|f\|^2+\|K_0-K \e^{\sigma}\|^2-\int_M f(K_0-K
\e^{\sigma}) \d\mu\\[.2pc]
&=\bigg\|-\frac{1}{2}f+K_0-K \e^{\sigma}\bigg\|^2.
\end{align*}

We will show that this is an equality since the inequality contradicts
the fact that $\{\sigma_n\}$ was a minimizing sequence. That is we can
construct a sequence $\{\tau \} \in C^{\infty}(M)$ such that $S[\tau]$
gets as close to $\big\|-\frac{1}{2}f+K_0-K\e^{\sigma}\big\|^2$ as we like.

Namely, for any $\epsilon>0$ we can construct, by the density of
$C^{\infty}$ in $W^{2,2}$, a function $\tau\in C^{\infty}(M)$
approximating $\sigma \in W^{2,2}$ such that
$\|\Delta_h\tau-f\|<\epsilon$ and $\|4(\e^{\tau}- \e^{\sigma})\|<
\epsilon/2$. Since
\begin{equation*}
S_{\tau} = \lim\limits_{n \rightarrow \infty}
S[\tau]=\bigg\|-\frac{1}{2}\Delta_h\tau+K_0 -K \e^{\tau}\bigg\|^2,
\end{equation*}
we have\vspace{.2pc}
\begin{equation*}
\bigg\| \sqrt{S_{\tau}}-\bigg \|-\frac{1}{2}f+K_0-K \e^{\sigma} \bigg\| \leq
\bigg\|\frac{1}{2}(f-\Delta_h\tau)-K(\e^{\tau}-
\e^{\sigma})\bigg\|\leq\epsilon.
\end{equation*}

$\left.\right.$\vspace{-1.1pc}

\noindent Now setting $\delta=\sqrt{S_0}-\|-\frac{1}{2}f+K_0-K
\e^{\sigma}\|>0$ and choosing $\epsilon <\delta/2$, and using
$\sqrt{S_{\tau}}\leq\|-\frac{1}{2}f+K_0-K \e^{\sigma}\|+ \epsilon$
we get, $\sqrt{S_{\tau}}<\sqrt{S_0}-\frac{\delta}{2}$ -- a
contradiction, since $S_0$ is the infimum of the functional.

Thus, $\lim_{n\rightarrow\infty}\|\Delta_h \sigma_n\|=\|f\|$, so
that, in fact, by Lemma~2.5, the convergence is in the strong
$L^2$ topology. This proves part (c).\hfill$\Box$
\end{proof}

\section{Smoothness and uniqueness}

Here we complete the proof of the main Theorem~2.1 by showing that\vspace{.2pc}

\setcounter{equation}{0}
\setcounter{defin}{0}
\begin{proposition}$\left.\right.$%\vspace{.5pc}
\begin{enumerate}
\renewcommand{\labelenumi}{\rm (\alph{enumi})}
\item The minimizing function $\sigma\in C^0(M)$ is smooth
and unique and corresponds to a metric of  negative curvature $K$.

\item $\sigma$  extends to $M - \cup_{i=1}^{n} p_i$ by zero.
\end{enumerate}\vspace{-1.5pc}
\end{proposition}

\pagebreak

\begin{proof}
(a) Let $b=(K_0-\frac{1}{2}\Delta_h^{\rm
distr}\sigma-K\e^{\sigma})\in L^2(M)$; note that $b=0$ on the boundary
$\partial M$, since $\sigma|_{\partial M} = 0$ and by redefining
$\Delta^{\rm distr}_h \sigma$ on the boundary, we can take $\Delta^{\rm
distr}_h \sigma|_{\partial M} = 0$ and we know that $K|_{\partial M} =
K_0$. According to Proposition~2.4, (c) $S_0= \int_M b^2\d\mu$.
Set $G(t)=S(\sigma+t\beta)-S_0$, where $\beta\in \{ f \in C^{\infty}(M)
| \partial_{\nu} f|_{\partial M} = 0, f|_{\partial M} =0 \}$. $G(t)$ for
fixed $\beta$ is smooth, $G(0)=0$ and $G(t)\geq 0$ for all $t$.
Therefore,
\begin{equation*}
\frac{\d G}{\d t}\bigg|_{t=0}=0.
\end{equation*}
A simple calculation yields
\begin{equation*}
\frac{\d G}{\d t}\bigg|_{t=0}=\int_M(-b\Delta_h\beta-2K\e^{\sigma}b\beta) \d\mu.
\end{equation*}
Now, $\int_M b \Delta_h \beta\ \d \mu - \int_M \beta \Delta_h b\ \d \mu
= \int_{\partial M} b \partial_{\nu} \beta\ \d l - \int_{\partial M}
\beta \partial_{\nu} b\ \d l =0$. Thus $b\in L^2(M)$ satisfies, in a
distributional sense, the following equation:
\begin{equation} \label{eq}
 -\Delta_h b-2K\e^{\sigma} b=0.
\end{equation}

First, we will show that $b=0$ is the only weak $L^2$ solution to
eq.~\eqref{eq}. Indeed, by elliptic regularity $b$ is smooth, so that multiplying
\eqref{eq} by $b$ and integrating over $M$ using the Stokes formula, we get
\begin{equation*}
\int_M \nabla b \cdot \nabla b -\int_{\partial M} \partial_{\nu} b \cdot
b - \int_M 2 K b^2\e^{\sigma}\d\mu=0.
\end{equation*}
The second term drops since we redefined $\Delta^{\rm distr} \sigma$ to be zero on the boundary and had $K=K_0$ on the boundary, thus   $b=0$ on $\partial M$. Recalling
that $K<0$ on $M$, all the remaining terms are positive or zero,
$b=0$ everywhere. Thus we have shown that $S_0=0$.

Secondly, equation $b=0$ for the minimizing function $\sigma\in
C^0(M)$ reads
\begin{equation} \label{boot}
\frac{1}{2}\Delta_h^{\rm distr}\sigma=K_{0}-K\e^{\sigma}\in C^{0}(M).
\end{equation}
Therefore, $\Delta_h^{\rm distr}\sigma$ belongs to $L^p(M)$ so that
$\sigma\in W^{2,p}$ for all $p$. By the Sobolev embedding theorem
it follows that $\sigma\in C^{1,\alpha}(M)$ for some
$0<\alpha<1$. Therefore, the right-hand side of eq.~\eqref{boot}
actually belongs to the space $C^{1,\alpha}(M)$, and
therefore $\sigma\in C^{3,\alpha}(M)$ and so on. This kind of
bootstrapping argument shows that $\sigma$ is smooth.

Equation $b \equiv 0$ satisfied by $\sigma$ now translates to
$K(\sigma) \equiv K$, where $K(\sigma)$ is the Gaussian curvature of the
metric ${\rm e}^{\sigma}h\ \d z \otimes \d \bar{z}$, $\sigma \in C^{\infty}(M)$.

Next to show uniqueness, let $\eta$ be another solution in class $T$. Then
$\eta$ satisfies
\begin{equation} \label{boot}
\frac{1}{2}\Delta_h^{\rm distr}\eta=K_{0}-K\e^{\eta}\in C^{0}(M)
\end{equation}
so that
\begin{equation*}
\Delta_h (\sigma - \eta) = - 2 K (\e^{\sigma} - \e^{\eta}).
\end{equation*}
Multiplying this equation by $\sigma - \eta$ and remembering that $\sigma$
and $\eta$ and their normal derivatives vanish on the boundary, we get
\begin{equation*}
- \int_M \d \zeta \wedge * \d \zeta = \int_M -2K (\sigma - \eta)(\e^{\sigma} -
\e^{\eta}) \d \mu,
\end{equation*}
where we set $\zeta = \sigma - \eta$. Since $-2K(\sigma - \eta) (\e^{\sigma} -
\e^{\eta}) \geq 0$, we conclude that $\d \zeta = 0$ and in fact $\zeta = 0$.

(b) $\sigma$ is in class $T$, i.e. $\sigma|_{\partial M} = \partial_{\nu}
\sigma |_{\partial M} = \Delta_h \sigma |_{\partial M} = 0$. Supposing we
minimize over the class of functions all of whose derivatives are zero at
$\partial M$. Let $\tau$ be the smooth function to which this minimizing
sequence converge to. By uniqueness of the solution, $\tau = \sigma$.
We can extend $\tau$ by zero on the discs $\cup_{i=1}^{n} D_i$.
Thus we can extend $\sigma$ by zero as well. \hfill $\Box$
\end{proof}

\section{An example of the initial metric}\vspace{.2pc}

Let $\Sigma - \{ p_i \}_{i=1}^{n}$ be the Riemann surface with
punctures. Let $D_i= \{ (r_i, \theta_i): 0< r_i < 1 \}$ be a disc about
$p_i$ such that $p_i$ is the only puncture on $D_i$. Recall $M = \Sigma
- \cup_{i=1}^{n} D_i $. $\partial M = \cup_{i=1}^{n} \partial D_i =
\cup_{i=1}^{n} \{ (r_i, \theta_i): r_i = 1 \}$. Let $D_i^{\prime} = \{
(r_i, \theta_i): 0 <r_i < 1.5 \} $ and $D_i^{\prime \prime} = \{ (r_i,
\theta_i): 0< r_i < 2 \}$. Let $f_i\ \d z \otimes \d \bar{z}$ be a
metric on $D_i^{\prime \prime}- \ p_i$. Let $\tilde{h}\ \d z \otimes \d
\bar{z}$ be a metric on $M$. We wish to construct a metric which will
interpolate between the metric on $M$ and the metric on $D_i^{\prime
\prime} $, $i=1,\ldots,n$. Let
\begin{align*}
\rho_i &= 1\quad \hbox{on}\  \{ 0 \leq r_i \leq 1.5 \}\\
       &= 0\quad \hbox{on}\  \{ r_i \geq 2 \}\\
&\quad\ \, 0 < \rho < 1 \quad \hbox{on}\ \{ 1.5 < r_i < 2 \}\\
\intertext{and}
\tilde{\rho} &= 1\quad \hbox{on}\ \Sigma - \cup_{i=1}^{n} D_i^{\prime \prime}\\
&= 0 \quad \hbox{on}\ D_i^{\prime}\ \forall i\\
& \quad\ \, 0 <  \tilde{\rho} < 1 \quad \hbox{on}\ \{ 1.5 < r_i <2 \}.
\end{align*}

Let
\begin{equation*}
h\ \d z \otimes \d \bar{z} = \sum\limits_{i=1}^{n} \rho_i f_i\ \d z \otimes \d \bar{z} + \tilde{\rho} \tilde{h}\ \d z \otimes \d \bar{z}.
\end{equation*}
This is a metric because $h$ is positive everywhere. On $D_i^{\prime}$,
$h\ \d z \otimes \d \bar{z} = f_i\ \d z \otimes \d \bar{z}$ and on
$\Sigma - \cup_{i=1}^{n} D_i^{\prime \prime}$, $h\ \d z \otimes d
\bar{z} = \tilde{h}\ \d z \otimes \d \bar{z}$.

We can choose $f_i$ such that the Gaussian curvature $K_0$ of $h\ \d z
\otimes \d \bar{z}$ is negative on $\partial M$. One such example would
be to choose $f_i\ \d z \otimes \d \bar{z} $ on the disc
$D_i^{\prime}-p_i $\ \,to be\ \,$ -\log(\frac{r_i}{4}) [(\d r_i)^2 + r_i^2 (\d
\theta_i)^2 ] $, where $(r_i, \theta_i )$ are standard polar coordinates
on the disc. This metric has a true singularity at the puncture. $K_0 =
\frac{1}{2 \sqrt{EG}}
\left[\left(\frac{E_{\theta_i}}{\sqrt{EG}}\right)_{\theta_i} +
\left(\frac{G_{r_i}}{\sqrt{EG}}\right)_{r_i}\right]$ where $E=
-\log\left(\frac{r_i}{4}\right)$, $G= -r_i^2
\log\left(\frac{r_i}{4}\right)$ \cite{DoC}. On $D_i$, $K_0 = \frac{1}{2
r_i^2 \left(\log (r_{i}/4)\right)^3}$ and therefore
negative. In particular, on $\partial M = \cup_{i=1}^{n} \{(r_i,
\theta_i): r_i=1 \}$, $K_0$ is negative.

\section*{Acknowledgement}\vspace{.2pc}

The author would like to thank V~Pati and E~K~Narayanan for very useful
discussions. The author would also like to thank the referee for his
comments.

\end{document}